\documentclass[aos,MSNbibl,seceqn,nameyear,dvips]{arximspdf}
\usepackage{algorithm}
\usepackage{dcolumn}
\usepackage{multirow}
\usepackage{graphicx}

\doi{10.1214/12-AOS1017} 
\volume{40}
\issue{3}
\pubyear{2012}
\firstpage{1403}
\lastpage{1429}

\makeatletter

\newcolumntype{d}[1]{D{.}{.}{#1}}

\newproclaim{remark}{Remark}

\newtheorem{theorem}{Theorem}[section]
\newtheorem{cor}{Corollary}[section]

\newcommand{\boldbeta}{{\bolds{\beta}}}
\newcommand{\boldalpha}{{\bolds{\alpha}}}
\newcommand{\boldzero}{{\bolds{0}}}
\newcommand{\boldx}{{\mathbf{x}}}
\newcommand{\pr}{\operatorname{Pr}}

\newcommand{\bA}{\mathbf{A}}
\newcommand{\bX}{\mathbf{X}}

\makeatother

\begin{document}
\begin{frontmatter}

\title{Nonconcave penalized composite conditional likelihood
estimation of sparse Ising models\thanksref{T1}}
\runtitle{Penalized estimation of sparse Ising models}

\thankstext{T1}{Supported in part by NSF Grants DMS-08-46068
and DMS-08-54970.}

\begin{aug}
\author[A]{\fnms{Lingzhou} \snm{Xue}\ead[label=e1]{lzxue@stat.umn.edu}},
\author[A]{\fnms{Hui} \snm{Zou}\corref{}\ead[label=e2]{zouxx019@umn.edu}}
\and
\author[B]{\fnms{Tianxi} \snm{Cai}\ead[label=e3]{tcai@hsph.harvard.edu}}
\runauthor{L. Xue, H. Zou and T. Cai}
\affiliation{University of Minnesota, University of Minnesota and
Harvard University}
\address[A]{L. Xue\\
H. Zou\\
School of Statistics\\
University of Minnesota\\
Minneapolis, Minnesota 55455\\
USA\\
\printead{e1}\\
\hphantom{E-mail: }\printead*{e2}}
\address[B]{T. Cai\\
Department of Biostatistics\\
Harvard University\\
Boston, Massachusetts 02115\\
USA\\
\printead{e3}} 
\end{aug}

\received{\smonth{9} \syear{2011}}
\revised{\smonth{5} \syear{2012}}

%
\begin{abstract}
The Ising model is a useful tool for studying complex interactions
within a system. The estimation of such a model, however, is rather
challenging, especially in the presence of high-dimensional parameters.
In this work, we propose efficient procedures for learning a~sparse
Ising model based on a penalized composite conditional likelihood with
nonconcave penalties. \mbox{Nonconcave} penalized likelihood estimation has
received a lot of attention in recent years. However, such an approach
is computationally prohibitive under high-dimensional Ising models. To
overcome such difficulties, we extend the methodology and theory of
nonconcave penalized likelihood to penalized composite conditional
likelihood estimation. The proposed method can be efficiently
implemented by taking advantage of coordinate-ascent and
minorization--maximization principles. Asymptotic oracle properties of
the proposed method are established with NP-dimensionality. Optimality
of the computed local solution is discussed. We demonstrate its finite
sample performance via simulation studies and further illustrate our
proposal by studying the Human Immunodeficiency Virus type 1 protease
structure based on data from the Stanford HIV drug resistance database.
Our statistical learning results match the known biological findings
very well, although no prior biological information is used in the data
analysis procedure.
\end{abstract}

%

\end{frontmatter}

\section{Introduction}\label{sec1}

The Ising model was first introduced in statistical phys\-ics
[\citet{Isi25}] as a mathematical model for describing magnetic
interactions and the structures of ferromagnetic substances. Although
rooted in physics, the Ising model has been successfully exploited to
simplify complex interactions for network exploration in various
research fields such as social-economics [\citet{Sta08}], protein
modeling [\citet{IrbPetPot96}] and statistical genetics
[\citet{MajLiOtt01}]. Following the terminology in physics,
consider an Ising model with $K$ magnetic dipoles denoted by $X_j$, $1
\leq j \leq K$. Each $X_j$ equals $+1$ or $-1$, corresponding to the up
or down spin state of the $j$th magnetic dipole. The energy function is
defined as $E=-\sum_{i \neq j}\beta_{ij}\frac{X_iX_j}{4}$, where the
coupling coefficient $\beta_{ij}$ describes the physical interactions
between dipoles $i$ and $j$ under the external magnetic field,
$\beta_{ii}=0$ and $\beta_{ij}=\beta_{ji}$ for any $(i,j)$. According
to Boltzmann's law, the joint distribution of $\bX= (X_1, \ldots, X_K)$
should be
%
%
\begin{equation}
\label{ising} \pr(X_{1}=x_{1}, \ldots, X_{K}=x_{K})=
\frac{1}{Z(\boldbeta)}\exp\biggl(\sum_{(i,j)}
\frac{\beta_{ij}
x_{j}x_{i}}{4}\biggr),
\end{equation}
where $Z(\boldbeta)$ is the partition function.

In this paper we focus on learning sparse Ising models; that is, many
coupling coefficients are zero. Our research is motivated by
the HIV drug resistance study where understanding the inter-residue
couplings (interactions)
could potentially shed light on the mechanisms of drug resistance. A
suitable statistical learning method
is to fit a sparse Ising model to the data, in order to discover the
inter-residue couplings.
More details are given in Section \ref{sec5}. In the recent
statistical literature, penalized likelihood estimation has become
a standard tool for sparse estimation. See a recent review paper
by \citet{FanLv10}. In principle we can follow the penalized
likelihood estimation paradigm to derive a sparse penalized
estimator of the Ising model. Unfortunately, the penalized
likelihood estimation method is very difficult to compute under
the Ising model because the partition function $Z(\boldbeta)$ is
computationally intractable when the number of dipoles is
relatively large. On the other hand, the composite likelihood idea
[\citet{Lin88}, \citet{VarReiFir11}] offers a nice
alternative. To elaborate, suppose we have $N$ independent
identically distributed (i.i.d.) realizations of $\bX$ from the Ising
model, denoted by $\{(x_{1n},\ldots, x_{Kn}), n = 1,\ldots, N\}$.
Let $\theta_{j}= P(X_i=x_j|\bX_{(-j)})$, describing
the conditional distribution of the $j$th dipole given the remaining
dipoles, where $\bX_{(-j)}$ denotes $\bX$ with the $j$th element removed.
By (\ref{ising}), it is easy see that for the $n$th
observation,\looseness=1
\[
\theta_{jn} = \frac{\exp(\sum_{k\dvtx k \neq j}
\beta_{jk}x_{jn}x_{kn})}{\exp(\sum_{k\dvtx k\neq j}
\beta_{jk}x_{jn}x_{kn})+1}.
\]\looseness=0
Note that $\theta_{jn}$ does not involve the partition function.
The conditional log-likelihood of the $j$th dipole, given the
remaining dipoles, is given by
\[
\ell^{(j)}=\frac{1}{N}\sum^N_{n=1}
\log(\theta_{jn}).
\]
As in \citet{Lin88} a composite log-likelihood function can
be defined as
\[
\ell_c = \sum^K_{j=1}
\ell^{(j)}.
\]
This kind of composite conditional likelihood was also called
pseudo-likeli\-hood in \citet{Bes74}. Another popular type of
composite likelihood is composite marginal likelihood
[\citet{Var08}]. Maximum composite likelihood is especially useful
when the full likelihood is intractable. Such an approach has important
applications in many areas including spatial statistics, clustered and
longitudinal data and time series models. A nice review on the recent
developments in composite likelihood can be found in
\citet{VarReiFir11}.

To estimate a high-dimensional sparse Ising
model, we consider the following penalized composite
likelihood estimator:
%
%
\begin{equation}
\label{PCL} \widehat\boldbeta=\mathop{\arg\max}_{\boldbeta} \Biggl\{\ell_c
(\boldbeta)-\sum^K_{j=1}\sum
^K_{k=j+1}P_{\lambda}\bigl(|\beta_{jk}|\bigr)
\Biggr\},
\end{equation}
where $P_{\lambda}(t)$ is a positive penalty function defined on
$[0,\infty)$. In this work we focus primarily on the LASSO penalty
[\citet{Tib96}] and smoothly clipped absolute deviation (SCAD)
penalty [\citet{FanLi01}]. The LASSO penalty is
$P_{\lambda}(t)=\lambda t$. The SCAD penalty is defined by
\[
P^{\prime}_{\lambda}(t)=\lambda\biggl\{I(t \le\lambda)+
\frac{(a\lambda-t)_{+}}{(a-1)\lambda}I(t>\lambda) \biggr\},\qquad t \ge0; a>2.
\]
Following \citet{FanLi01} we set $a=3.7$. We should make it clear
that when $P_{\lambda}(t)$ is nonconcave, $\widehat\boldbeta$
should be
understood as a good local maximizer of (\ref{PCL}). See discussions
in Section \ref{sec2}.

The optimization problem in (\ref{PCL}) is very challenging
because of two major issues: (1) the number of unknown parameters
is $\frac{1}{2}K(K-1)$, and hence the optimization problem is
high dimensional in nature; and (2) the penalty function is
concave and nondifferentiable at zero, although $\ell_c$ is a
smooth concave function. We propose to combine the strengths of
coordinate-ascent and minorization--maximization, which
results in two new algorithms, CMA and LLA--CMA, for computing a local
solution of the nonconcave penalized composite likelihood.
See Section \ref{sec2} for details. With the aid of the new algorithms,
the SCAD
penalized estimators are able to enjoy computational efficiency
comparable to that of the LASSO penalized estimator.

\citet{FanLi01} advocated the oracle properties of the nonconcave
penalized likelihood estimator in the sense that it performs as well as
the oracle
estimator which is the hypothetical maximum likelihood estimator
knowing the true
submodel. Zhang (\citeyear{Zha10N1}) and \citet{LvFan09} were among
the first to study the concave penalized
least-squares estimator with NP-dimensionality ($p$ can grow faster
than any polynomial function of $n$).
\citet{FanLv11} studied the asymptotic properties of nonconcave penalized
likelihood for generalized linear models with NP-dimensionality.
In this paper we show that the oracle model selection theory remains to
hold nicely for
nonconcave penalized composite likelihood with NP-dimensionality.
Furthermore, we show that under certain regularity conditions the
oracle estimator can be attained asymptotically via the LLA--CMA algorithm.

There is some related work in the literature. \citet{RavWaiLaf10}
viewed the Ising model as a binary
Markov graph and used a neighborhood LASSO-penalized logistic
regression algorithm to select the edges. Their idea is an
extension of neighborhood selection by LASSO regression proposed
by \citet{MeiBuh06} for estimating Gaussian
graphical models. \citet{HofTib09} suggested
using the LASSO-penalized pseudo-likelihood to estimate binary
Markov graphs. However, they did not provide any
theoretical result nor application. In this paper we compare the
LASSO and the SCAD penalized composite likelihood estimators and
show the latter has substantial advantages with respect to both
numerical and theoretical properties.

The rest of this paper is organized as follows. In Section \ref{sec2}, we
introduce the CMA and LLA--CMA algorithms. The
statistical theory is presented in Section \ref{sec3}. Monte Carlo
simulation results are shown in Section~\ref{sec4}. In Section~\ref
{sec5} we present a
real application of the proposed method to study the network
structure of the amino-acid sequences of retroviral proteases
using data from the Stanford HIV drug resistance database.
Technical proofs are relegated to the \hyperref[app]{Appendix}.

\section{Computing algorithms}\label{sec2}
In this section we discuss how to efficiently implement the penalized
composite likelihood estimators. As mentioned before, the computational
challenges come from (1) penalizing the concave composite likelihood
with a nonconcave penalty which is not differentiable at zero; (2) the
intrinsically high dimension of the unknown parameters.
\citet{ZouLi08} proposed the local linear approximation (LLA)
algorithm to derive an iterative \mbox{$\ell_1$-optimization} procedure for
computing nonconcave penalized estimators. The basic idea behind LLA is
the minorization--maximiza\-tion principle [\citet{LanHunYan00},
\citet{HunLan04}, \citet{HunLi05}]. Coordinate-ascent (or descent) algorithms
[\citet{Tse88}] have been successfully used for solving penalized
estimators with LASSO-type penalties; see, for example,
\citet{Fu98}, \citet{DauDefDeM04},
\citet{GenLewMad07}, \citet{YuaLin06}, \citet{MeivanBuh08},
\citet{WuLan08} and \citet{FriHasTib10}. In this paper we
combine the strengths of
minorization--maximization and coordinatewise optimization to overcome
the computational challenges.

\subsection{The CMA algorithm}\label{sec21}

Let $\widetilde\boldbeta$ be the current estimate. The
coordinate-ascent algorithm sequentially updates $\widetilde
\beta_{ij}$ by solving the following univariate optimization
problem:
%
%
\begin{equation}
\label{sec2eq1}
\qquad \widetilde\beta_{jk} \Leftarrow\mathop{\arg\max}_{\beta_{jk}}
\bigl\{\ell_c\bigl(\beta_{jk};\beta_{j'k'}=\widetilde
\beta_{j'k'}, \bigl(j',k'\bigr) \neq(j,k)
\bigr)-P_{\lambda}\bigl(|\beta_{jk}|\bigr) \bigr\}.
\end{equation}
However, we do not have a closed-form solution for the maximizer
of (\ref{sec2eq1}). The exact maximization has to be conducted by
some numerical optimization routine, which may not be a
good choice in the coordinate-ascent algorithm \mbox{because} the
maximization routine needs to be repeated many times to reach
convergence. On the other hand, one can find an update to increase,
rather than maximize, the objective function in (\ref{sec2eq1}),
maintaining the crucial ascent property of the coordinate-ascent
algorithm. This idea is in line with the generalized EM algorithm
[\citet{DemLaiRub77}] in which
one seeks to increase the expected log likelihood in the M-step.

First, we observe that for
any $\beta_{ij}$
%
%
\begin{equation}
\label{bd1} \frac{\partial^2 \ell_c(\boldbeta)}{\partial
\beta_{jk}^2}=-\frac{1}{N}\sum
_{n=1}^{N} \bigl(\theta_{kn}(1-
\theta_{kn})+ \theta_{jn}(1-\theta_{jn})\bigr)\ge-
\frac{1}{2}.
\end{equation}
Thus, by Taylor's expansion, we have
\[
\ell_c\bigl(\beta_{jk};\beta_{j'k'}=\widetilde
\beta_{j'k'}, \bigl(j',k'\bigr) \neq(j,k)\bigr)
\geq Q(\beta_{jk}),
\]
where
%
%
\begin{eqnarray}
\label{Q} Q(\beta_{jk}) &\equiv& \ell_c\bigl(
\beta_{jk}=\widetilde\beta_{jk};\beta_{j'k'}=\widetilde
\beta_{j'k'}, \bigl(j',k'\bigr) \neq(j,k)\bigr)
\nonumber\\[-8pt]\\[-8pt]
&&{} +\widetilde z_{jk}(\beta_{jk}-\widetilde{\beta}_{jk})-
\tfrac
{1}{4}(\beta_{jk}-\widetilde{\beta}_{jk})^2,\nonumber
\\
%
%
\widetilde z_{jk}&=&\frac{\partial\ell_c(\boldbeta)}{\partial
\beta_{jk}}\bigg\vert_{\boldbeta=\widetilde{\boldbeta}}=
\frac{1}{N} \sum_{n=1}^{N}
x_{kn}x_{jn}\bigl(2-\theta_{kn}(\widetilde{\beta})
- \theta_{jn}(\widetilde{\beta})\bigr).
\end{eqnarray}
Next, \citet{ZouLi08} showed that
%
%
\begin{equation}
\label{L} P_{\lambda}\bigl(|\beta_{jk}|\bigr) \le P_{\lambda}\bigl(|\widetilde
\beta_{jk}|\bigr)+ P^{\prime}_{\lambda}\bigl(|\widetilde
\beta_{jk}|\bigr)\cdot\bigl(|\beta_{jk}|-|\widetilde\beta_{jk}|\bigr)
\equiv L\bigl(|\beta_{jk}|\bigr).
\end{equation}
Combining (\ref{Q})--(\ref{L}) we see that $Q(\beta_{jk})-
L(|\beta_{jk}|)$ is a minorization function of the objective
function in (\ref{sec2eq1}). We update
$\widetilde\beta_{jk}$ by
%
%
\begin{equation}
\label{sec2eq2} \widetilde\beta_{jk}^{\mathrm{new}} = \mathop{\arg
\max}_{\beta_{jk}} \bigl\{Q(\beta_{jk})-L\bigl(|\beta_{jk}|\bigr) \bigr
\},
\end{equation}
whose solution is given by
$ \widetilde\beta_{jk}^{\mathrm{new}} = S(\widetilde
\beta_{jk}+2\widetilde z_{jk},2 P'_{\lambda}(|\widetilde\beta_{jk}|) )
$
where $S(r,t)=\operatorname{sgn}(r)(|r|-t)_{+}$ denotes the
soft-thresholding operator
[\citet{Tib96}].
The above arguments lead to Algorithm \ref{algo1} below, which we call the
coordinate-minorization-ascent (CMA) algorithm.

\begin{algorithm}[t]
\caption{The CMA algorithm}\label{algo1}
\begin{longlist}[(3)]\vspace*{3pt}
\item[(1)] Initialization of $\widetilde\boldbeta$.

\item[(2)] Cyclic coordinate-minorization-ascent: sequentially update
$\widetilde\beta_{ij}$ ($1 \leq j <k \leq K$) via soft-thresholding $
\widetilde\beta_{jk} \Leftarrow S(\widetilde\beta_{jk}+2\widetilde z_{jk},2
P'_{\lambda}(|\widetilde\beta_{jk}|) ). $

\item[(3)] Repeat the above cycle till convergence.
\end{longlist}
\end{algorithm}

%
\begin{remark}\label{remark1}
It is easy to prove that Algorithm \ref{algo1} has a
nice ascent property which is a direct consequence of the
minorization--maximizaton principle. Note that Algorithm \ref{algo1} can be
directly used to compute the LASSO-penalized composite likelihood
estimator. We simply
modify the coordinate-wise updating formula as $
\widetilde\beta_{jk} \Leftarrow S(\widetilde\beta_{jk}+2\widetilde
z_{jk},2\lambda)$.

In practice we need to specify the $\lambda$ value.
BIC has been shown to perform very well for selecting the tuning
parameter of the penalized likelihood estimator [\citet
{WanLiTsa07}]. The BIC score is defined as
%
%
\begin{equation}
\widehat\lambda=\mathop{\arg\max}_{\lambda} \biggl\{ 2\ell_c\bigl(\widehat
\boldbeta(\lambda)\bigr)-\log(n) \cdot\sum_{(j,k)}I
\bigl(\widehat\beta_{jk}(\lambda) \neq0\bigr)\biggr\}.
\end{equation}
BIC is used to tune all methods considered in this work. We use SCAD1
to denote the SCAD solution computed by Algorithm \ref{algo1} with the
BIC tuned
LASSO solution being the starting value.

For computational efficiency considerations, we implement Algorithm \ref{algo1}
by using the path-following idea and some other tricks, including
warm-starts and active-set-cycling [\citet{FriHasTib10}].
We have implemented the algorithm in R language functions. The core
cyclic coordinate-wise soft-thresholding operations were carried out in
C.
\end{remark}
%
%
\begin{remark}\label{remark2}
As suggested by a referee, the
coordinate-gradient-ascent (CGA) algorithm is a natural alternative to
Algorithm \ref{algo1} for solving the LASSO-penalized composite likelihood
estimator. The CGA algorithm has successfully used to solve other
penalized models. See \citet{GenLewMad07}, \citet{MeivanBuh08},
St{\"a}dler, B{\"u}hl\-mann and van~de Geer (\citeyear{StaBuhvan10}) and \citet{SchBuhvan11}.
In the CGA algorithm we need to find a good step size along the
gradient direction to guarantee the ascent property after each
coordinate-wise update. These extra computations are necessary
for the CGA algorithm, but are not needed in the CMA algorithm. We have
also implemented the CGA algorithm to solve the LASSO estimator and
found that the CMA algorithm is about five times faster than the CGA
algorithm. See Section \ref{sec4} for the timing comparison details.
\end{remark}


\subsection{Issues of local solution and the LLA--CMA algorithm}\label{sec22}

The objective function in (\ref{PCL}) is generally nonconcave if
a\vadjust{\goodbreak}
nonconcave penalty function is used.
Using Algorithm \ref{algo1} we find a local solution to (\ref{PCL}),
but there
is no guarantee that it is the global solution. A similar case is
\citet{SchBuhvan11} where the objective
function is the LASSO-penalized maximum likelihood of a
high-dimensional linear mixed-effects model, and the authors derived a
coordinate-wise gradient descent algorithm to find a~local solution.

It should not be considered as a special weakness of Algorithm \ref
{algo1} or
other coordinate-wise descent algorithm as
in \citet{SchBuhvan11} that the algorithm can only find a local
solution, because in the current literature there is no algorithm that
can guarantee to find the global solution of nonconcave
maximization (or nonconvex minimization) problems, especially when the
dimension is huge.
Consider, for example, the EM algorithm, which is perhaps the most
famous algorithm in statistical literature.
The EM algorithm often offers an elegant way
to fit some statistical models that are formulated as nonconcave
maximization problems. However, the EM algorithm provides a local
solution in general. A recent application of the EM algorithm to
high-dimensional modeling can be found in
\citet{StaBuhvan10} who considered a LASSO-penalized maximum likelihood
estimator of a high-dimensional linear regression model with
inhomogeneous errors that are modeled by a finite mixture of Gaussians.
To handle the computational challenges in their problem,
\citet{StaBuhvan10} proposed a generalized EM algorithm in which
a~coordinate
descent loop is used in the M-step and showed that the obtained
solution is a local solution.

Our numerical results show that in the penalized composite likelihood
estimation problem the SCAD performs much better than
the LASSO. To offer theoretical understanding of their differences, it
is important to show that the obtained local solution of the
SCAD-penalized likelihood
has better theoretical properties than the LASSO estimator. In Section
\ref{sec3} we establish the asymptotic properties of the LASSO
estimator and a
local solution of~(\ref{PCL}) with the SCAD penalty.
However, a general technical difficulty in nonconcave maximization
problems is to show that the computed local solution is the one local
solution with proven theoretical properties.
In \citet{StaBuhvan10} and \citet{SchBuhvan11}, nice
asymptotic properties are established for their proposed methods
but it is not clear whether the computed local solutions could have
those theoretical properties. The same issue exists in \citet{FanLv11}.

To circumvent the technical difficulty, we can consider combining the
LLA idea [\citet{ZouLi08}] and Algorithm \ref{algo1} to solve (\ref
{PCL}) with
a nonconcave penalty.
The LLA algorithm turns a nonconcave penalization problem into a
sequence of weighted LASSO penalization
problems. Similar ideas of iterative LLA convex relaxation have been
used in \citet{CanWakBoy08},\vadjust{\goodbreak} \citet{Zha10N2} and \citet
{BraFanWan11}.
Applying the LLA algorithm to (\ref{PCL}), we need to iteratively solve
%
%
\begin{equation}
\label{PCL2} \widehat\boldbeta{}^{(m+1)}=\mathop{\arg\max}_{\boldbeta}
\Biggl\{
\ell_c (\boldbeta)-\sum^K_{j=1}
\sum^K_{k=j+1}w_{jk}\cdot|
\beta_{jk}| \Biggr\}
\end{equation}
for $m=0,1,2,\ldots$ where $w_{jk}=P^{\prime}_{\lambda}(|\widetilde\beta
^{(m)}_{jk}|)$. Note that
Algorithm \ref{algo1} can be used to solve (\ref
{PCL2}) by simply modifying the coordinate-wise updating
formula as $
\widetilde\beta_{jk} \Leftarrow S(\widetilde\beta_{jk}+2\widetilde z_{jk},2
w_{jk}). $ Therefore, we have the following LLA--CMA algorithm for
computing a local solution of
(\ref{PCL}).

\begin{algorithm}[t]
\caption{The LLA--CMA algorithm}\label{algo2}
\begin{enumerate}[(2)]\vspace*{3pt}
\item[(1)] Initialize $\widetilde\boldbeta{}^{(0)}$, and compute
$w_{jk}=P^{\prime}_{\lambda}(|\widetilde\beta{}^{(0)}_{jk}|)$.
\item[(2)] For $m=0,1,2,3,\ldots\,$, repeat the LLA iteration:
\begin{enumerate}[(2.b)]
\item[(2.a)] Use Algorithm \ref{algo1} to solve $\widehat\boldbeta{}^{(m+1)}$
defined in (\ref{PCL2});
\item[(2.b)] Update the weights $w_{jk}$ by $P^{\prime}_{\lambda
}(|\widetilde
\beta{}^{(m+1)}_{jk}|)$.\vspace*{-2pt}
\end{enumerate}
\end{enumerate}
\end{algorithm}

In Section \ref{sec3} we show that if the LASSO estimator is $\widetilde
\boldbeta{}^{(0)}$, then under certain regularity conditions the LLA--CMA
algorithm finds the oracle estimator with high probability. These
results suggest that we should take the following steps to compute the
SCAD solution by the LLA--CMA algorithm.

\begin{center}
\textit{The proposed LLA--CMA procedure for computing a SCAD
estimator}:
\end{center}
\begin{longlist}[Step 3.]
\item[Step 1.] Use Algorithm \ref{algo1} to compute the LASSO solution
path and
find the LASSO estimator by BIC.
\item[Step 2.] Use the LASSO estimator as $\widetilde\boldbeta{}^{(0)}$
in the LLA--CMA algorithm to compute the solution path of the
first iteration and use BIC to tune the first step solution. Then use
the tuned first step solution as $\widetilde\boldbeta{}^{(0)}$ in the
LLA--CMA algorithm to compute the solution path and use BIC to select
$\lambda$. The resulting estimator is denoted by SCAD2.
\item[Step 3.] For the chosen $\lambda$ of SCAD2, use Algorithm \ref
{algo2} to
compute the fully converged SCAD solution with SCAD2 being the starting
value. Denote this SCAD solution by SCAD2$^{**}$.
\end{longlist}

The construction of SCAD2 follows an idea in \citet{BuhMei08}.
Based on our experience, SCAD2$^{**}$ works slightly better than SCAD2,
but the two are generally very close. Generally we recommend using
SCAD2$^{**}$ in real applications.

\section{Theoretical results}\label{sec3}

In this section we establish the statistical theory for
the penalized composite conditional likelihood estimator using the SCAD
and the LASSO penalty, respectively.
Such results allow us to compare the SCAD and the LASSO estimators
theoretically.\vadjust{\goodbreak}

In order to present the theory we need some necessary notation. {For
a~matrix $\bA=(a_{ij})$, we define the following matrix norms: the
Frobenius norm $\|\bA\|_F=\sqrt{\sum_{i,j}a_{ij}^2}$, the entry-wise
$\ell_{\infty}$ norm $\|\bA\|_{\max}=\max_{i,j}|a_{ij}|$ and the matrix
$\ell_\infty$ norm $\|\bA\|_{\infty}=\max_i\sum_j|a_{ij}|$.} Let
$\boldbeta^*=\{\beta_{jk}^*\dvtx j<k\}$ denote the\vspace*{1pt} true
coefficients, $\mathcal{A}=\{(j,k)\dvtx\beta_{jk}^*\neq0, j<k\}$ and
$s=|{\mathcal A}|$. Define\vspace*{2pt} $\rho(s,N)=\min_{(j,k) \in
\mathcal{A}}|\beta^*_{jk}|$ which represents the weakness of the
signal. Let $H$ be the Hessian matrix of $\ell_c$ such that
\[
H_{(j_1k_1),(j_2k_2)}=-\frac{\partial^2 \ell_c(\boldbeta)}{\partial
\beta_{j_1k_1}\,\partial\beta_{j_2k_2}},
\]
$1 \le j_1<k_1 \le K$ and
$1 \le j_2<k_2 \le K$. For simplicity we use $H^*=H(\boldbeta^*)$.
We partition $H$ and $\boldbeta$ according to $\mathcal A$ as
$
({H_{\mathcal{AA}} \atop H_{{\mathcal A}^c \mathcal{A}}}
\enskip{H_{\mathcal{A} {\mathcal A}^c} \atop H_{{\mathcal A}^c{\mathcal
A}^c}} )
$ and $\boldbeta= (\boldbeta_{{\mathcal A}}^T,
\boldbeta_{{\mathcal A}^c}^T)^T$, respectively. We
let
\[
\bX_{\mathcal A}=\bigl(X_{j}\dvtx(j,k) \mbox{ or }(k,j) \in
\mathcal{A}
\mbox{ for some } k\bigr)
\]
and
\[
\boldx_{\mathcal A n}=\bigl(x_{jn}\dvtx(j,k) \mbox{ or }(k,j)
\in\mathcal
{A} \mbox{ for some } k\bigr).
\]
Finally, we define
\begin{eqnarray*}
b&=&\lambda_{\min}\bigl(E\bigl[H^*_{\mathcal{A}\mathcal{A}}\bigr]\bigr),
\\
B&=&\lambda_{\max}\bigl(E\bigl[\bX_{\mathcal A}\bX^T_{\mathcal A}
\bigr]\bigr),
\\
\phi&=&\bigl\Vert E\bigl[H^{*}_{\mathcal{A}^c \mathcal
A}\bigr]\bigl(E
\bigl[H^{*}_{\mathcal{A} \mathcal{A}}\bigr]\bigr)^{-1}
\bigr\Vert_{\infty}.
\end{eqnarray*}

Define the oracle estimator as $\widehat{\boldbeta
}{}^{\mathrm{oracle}}=(\widetilde{\boldbeta}{}^{\mathrm{hmle}}_{\mathcal{A}},0)$
where
\[
\widetilde{\boldbeta}{}^{\mathrm{hmle}}_{\mathcal{A}}=
\mathop{\arg\max}_{\boldbeta_{\mathcal{A}}}
\ell_c\bigl((\boldbeta_{\mathcal{A}},0)\bigr).
\]
If we knew the true submodel, then we would use the oracle estimator to
estimate the Ising model.

%
%
\begin{theorem}\label{thm1}
Consider the SCAD-penalized composite likelihood defined in (\ref
{PCL}). We have the following two conclusions:
\begin{longlist}[(2)]
\item[(1)]
For any $R<\frac{b}{3B}\frac{\sqrt{N}}{s}$, we have
%
%
\begin{equation}
\label{scad1} \Pr\Biggl(\bigl\Vert\widetilde\boldbeta{}^{\mathrm{hmle}}_{\mathcal
A}-
\boldbeta^*_{\mathcal A} \bigr\Vert_2 \leq\sqrt{\frac{s}{N}} R
\Biggr) \ge1-\tau_1
\end{equation}
with $
\tau_1=\exp(-R^2\frac{b^2}{8^3})+2s^2\exp(-\frac{N}{s^2}\frac
{b^2}{2})+2s^2\exp(-\frac{N}{s^2}\frac{B^2}{8}).
$
\item[(2)] Pick a $\lambda$ satisfying
{$\lambda<\min(\frac{\rho(s,N)}{2a},\frac{{ (2\phi+1)} b^2}{3sB})$}.
With probability at
least $1-\tau_2$, $\widehat{\boldbeta}{}^{\mathrm{oracle}}$ is a local
maximizer\vadjust{\goodbreak} of the SCAD-penalized composite likelihood estimator where
%
%
\begin{eqnarray}
\label{tau2} \qquad\tau_2&=&\exp\biggl(-R^{2}_*
\frac{b^2}{8^3}\biggr) +K^2\exp\biggl(-\frac{N\lambda^2}{32{ (2\phi
+1)}^2}\biggr)
\nonumber
\\
&&{} +\exp\biggl(-\frac{N\lambda}{{ 3B(2\phi+1)}s}\frac{b^2}{8^3}\biggr) +K^2s
\exp\biggl(-\frac{Nb^2}{2s^3}\biggr)+2s^2\exp\biggl(-
\frac{b^2N}{8s^3}\biggr)
\\
&&{} +4s^2\biggl[\exp\biggl(-\frac{N}{s^2}\frac{b^2}{2}
\biggr)+\exp\biggl(-\frac{N}{s^2}\frac
{B^2}{8}\biggr)\biggr]\nonumber
\end{eqnarray}
and
$R_*=\min(\frac{1}{2}\sqrt{\frac{N}{s}}\rho(s,N),\frac
{b}{3B}\frac{\sqrt{N}}{s})$.
\end{longlist}
\end{theorem}

We also analyzed the theoretical properties of the LASSO estimator. If
the LASSO can consistently select the true model, it must equal to the
hypothetical LASSO estimator $(\widetilde{\boldbeta}_{\mathcal
{A}},0)$ where
\[
\widetilde{\boldbeta}_{\mathcal{A}}=\mathop{\arg\max}_{\boldbeta
_{\mathcal{A}}} \biggl\{
\ell_c\bigl((\boldbeta_{\mathcal{A}},0)\bigr)-\lambda\sum
_{(j,k)\in\mathcal{A}}|\beta_{jk}|\biggr\}.
\]

%
%
\begin{theorem}\label{thm2} Consider the LASSO-penalized composite
likelihood estimator.
\begin{longlist}[(2)]
\item[(1)] Choose $\lambda$ such that $\lambda s<\frac{8b^2}{3B}$.
$\Pr(\Vert\widetilde\boldbeta_{\mathcal
A}-\boldbeta^*_{\mathcal A} \Vert_2 \leq\frac{16\lambda
\sqrt{s}}{b} ) \ge1-\tau'_1$
with
\[
\tau'_1=e^{-N\lambda^2/2}+2s^2\biggl[\exp
\biggl(\frac{-Nb^2}{2s^2}\biggr)+\exp\biggl(\frac
{-NB^2}{8s^2}\biggr)\biggr].
\]
\item[(2)]Assume the ir-representable condition
$
\phi\leq1-\eta<1.
$ Choose $\lambda$
such that {$\lambda s
<\min(\frac{b^2}{16^2B}\frac{\eta/3}{{ 4-\eta}},\frac
{8b^2}{3B})$}. Then $(\widetilde\boldbeta_{{\mathcal A}},0)$ is the
LASSO-penalized composite likelihood estimator with probability at
least $1-\tau'_2$,
where
\begin{eqnarray*}
\tau'_2 &=& e^{-N\lambda^2/2}+K^2s\exp
\biggl(-\frac{Nb^2\eta^2}{8s^3}\biggr)+K^2\exp\biggl(- \frac{N\lambda
^2\eta^2}{32({ 4-\eta})^2}
\biggr)
\nonumber
\\
&&{} +2s^2\biggl[\exp\biggl(-\frac{Nb^2\eta^2}{2s^3(2-\eta)^2}\biggr)+\exp
\biggl(
\frac
{-Nb^2}{2s^2}\biggr)+\exp\biggl(\frac{-NB^2}{8s^2}\biggr)\biggr].
\end{eqnarray*}
\end{longlist}
\end{theorem}

In Theorems \ref{thm1} and \ref{thm2} the
three quantities $b$, $B$ and $\phi$ do not need to be constants.
We can obtain a more straightforward understanding of the
properties of the penalized composite likelihood estimators by
considering the asymptotic consequences of these probability
bounds. To highlight the main point, we consider $b$, $B$ and
$\phi$ are fixed constants and derive the following asymptotic
results.

%
%
\begin{cor}\label{thm3}
Suppose that $b$, $B$ and $\phi$ are fixed constants and further
assume $N \gg s^3\log(K)$ and $\rho(s,N) \gg
\sqrt{\frac{\log(K)}{{N}}}$.\vadjust{\goodbreak}
\begin{longlist}[(2)]
\item[(1)]
Pick the SCAD penalty parameter $\lambda^{\mathrm{scad}}$ satisfying
\[
\lambda^{\mathrm{scad}}<\min\biggl(\frac{\rho(s,N)}{2a},\frac{{ (2\phi+1)}
b^2}{3sB}\biggr),\qquad
\lambda^{\mathrm{scad}} \gg\sqrt{\frac{\log(K)}{{N}}}.
\]
With
probability tending to 1, the oracle estimator is a local
maximizer of the SCAD-penalized estimator and $\Vert\widehat
\boldbeta{}^{\mathrm{oracle}}_{\mathcal A}-\boldbeta^*_{\mathcal A} \Vert_2
=O_P(\sqrt{\frac{s}{N}})$.
\item[(2)] Assume the ir-representable condition in Theorem \ref
{thm2}. Pick the LASSO penalty parameter $\lambda^{\mathrm{lasso}}$ satisfying
\[
\min\biggl(\frac{1}{\sqrt{s}}\rho(s,N),\frac{1}{s}\biggr) \gg
\lambda^{\mathrm{lasso}} \gg\frac{1}{\sqrt{N}};
\]
then the LASSO estimator consistently selects the true model
and $\Vert\widehat\boldbeta{}^{\mathrm{lasso}}_{\mathcal
A}-\boldbeta^{*}_{\mathcal A}\Vert_2={O_P}(\lambda^{\mathrm{lasso}} \sqrt{s})$.
\end{longlist}
\end{cor}
%
%
\begin{remark}\label{remark3}
For the LASSO-penalized least squares, it is now known that the model
selection consistency critically depends on the ir-representable
condition [\citet{ZhaYu06}, \citet{MeiBuh06},
\citet{Zou06}]. A similar condition is again needed in the
LASSO-penalized composite likelihood. Furthermore, Corollary \ref{thm3}
shows that even when it is possible for the LASSO to achieve consistent
selection, $\lambda^{\mathrm{lasso}}$ should be much greater than
$\sqrt{\frac{1}{N}}$, which means that $\lambda^{\mathrm{lasso}} \sqrt{s}
\gg\sqrt{\frac{s}{N}}$. So the LASSO yields larger bias than the SCAD.
\end{remark}
%
%
\begin{remark}\label{remark4}
We have shown that asymptotically speaking
the oracle estimator is in fact a local solution of the
SCAD-penalized composite likelihood model. This property is stronger
than the oracle properties defined in \citet{FanLi01}.
Our result is the first to show that the oracle model selection theory
holds nicely for nonconcave
penalized composite conditional likelihood models with
NP-dimensionality. The usual composite likelihood theory in the literature
is only applied to the fixed-dimension setting. Our result fills a
long-standing gap in the composite likelihood literature.

What we have shown so far is the existence of a SCAD-penalized
estimator that is superior to the LASSO-penalized estimator.
Moreover, we would like to show that
the computed SCAD estimator is equal to the oracle estimator. As
discussed earlier in Section \ref{sec22}, such a result is very
difficult to
prove due to the nonconcavity
of the penalized likelihood function. See also \citet{FanLv11},
\citet{StaBuhvan10} and \citet{SchBuhvan11}.

If one can prove that the objective function has only one maximizer,
then the computed solution and the theoretically proven solution must
be the same.
This idea has been used in \citet{FanLv11} to study the nonconcave
penalized generalized linear models and \citet{BraFanJia} to study
the nonconcave penalized Cox proportional hazards models.
Their arguments are based on the observation that the SCAD penalty
function has a finite maximum concavity [\citet{Zha10N1},
\citet{LvFan09}]. Hence, if the smallest eigenvalue of the Hessian
matrix of the
negative log-likelihood is sufficiently large, the overall penalized
likelihood function is concave and hence has a unique global maximizer.
This argument requires that the sample size is greater than the
dimension; otherwise,
the Hessian matrix does not have full rank. To deal with the
high-dimensional case, \citet{FanLv11} further refined their
arguments by considering a subspace denoted by $\mathbb{S}_s$, which
is the union of
all $s$-dimensional coordinate subspaces. Under some regularity
conditions, \citet{FanLv11} showed that the oracle estimator is the
unique global maximizer in $\mathbb{S}_s$,
which was referred to as restricted global optimality. Then by assuming
that the computed solution has exactly $s$ nonzero elements, it can be
concluded that the computed solution is in
$\mathbb{S}_s$ and hence equals the oracle estimator; see Proposition
3.b of \citet{FanLv11}. However, a fundamental problem with these
arguments is that we have no idea whether the computed solution selects
$s$ nonzero coefficients, because $s$ is unknown.

Here we take a different route to tackle the local solution issue.
Instead of trying to prove the uniqueness of maximizer,
we directly analyze the local solution by the LLA--CMA algorithm and
discuss under which regularity conditions the LLA--CMA algorithm can
actually find the oracle estimator.
\end{remark}

%
%
\begin{theorem}\label{thm4}
Consider the SCAD-penalized composite likelihood estimator in (\ref
{PCL}). Let $\widehat\boldbeta{}^{\mathrm{scad}}$ be the local solution
computed by Algorithm \ref{algo2} (the LLA--CMA algorithm) with
$\widetilde
\boldbeta{}^{(0)}$ being the initial value.
Pick a $\lambda$ satisfying {$\lambda<\min(\frac{\rho
(s,N)}{2a},\frac{{ (2\phi+1)} b^2}{3sB})$}.
Write $\tau_0=\Pr(\Vert\widetilde\boldbeta{}^{(0)}-\boldbeta^*
\Vert_{\infty}>\lambda)$.
\begin{longlist}[(2)]
\item[(1)] The LLA--CMA algorithm finds the oracle estimator after one
LLA iteration with probability at least $1-\tau_0-\tau_3$ where
\begin{eqnarray*}
\tau_3 &=& K^2\exp\biggl(\frac{-N\lambda^2}{32{ (2\phi+1)}^2}\biggr
)+{\exp
\biggl(\frac
{-N\lambda}{{ 3B(2\phi+1)} s}\frac{b^2}{8^3}\biggr)}+K^2s\exp\biggl(
\frac
{-Nb^2}{2s^3}\biggr)
\nonumber
\\
&&{} +2s^2\biggl[\exp\biggl(-\frac{Nb^2}{8s^3}\biggr)+\exp\biggl(-
\frac{N}{s^2}\frac
{b^2}{2}\biggr)+\exp\biggl(-\frac{N}{s^2}
\frac{B^2}{8}\biggr)\biggr].
\end{eqnarray*}

\item[(2)] The LLA--CMA algorithm converges after two LLA iterations
and~$\widehat\boldbeta{}^{\mathrm{scad}}$ equals the oracle estimator with
probability at least
$1-\tau_0-\tau_2$,
where~$\tau_2$ is defined in (\ref{tau2}).
\end{longlist}
\end{theorem}

Theorem \ref{thm4} can be used to drive the following asymptotic result.
%
%
\begin{cor}\label{thm5}
Suppose that $b$, $B$ and $\phi$ are fixed constants, and further
assume $N \gg s^3\log(K)$ and $\rho(s,N) \gg
\frac{\max(\sqrt{\log(K)},16\sqrt{s}/b)}{\sqrt{N}}$.\vadjust{\goodbreak} Consider
the SCAD-penalized
composite likelihood estimator with the SCAD penalty parameter $\lambda
^{\mathrm{scad}}$ satisfying
\[
\lambda^{\mathrm{scad}}<\min\biggl(\frac{\rho(s,N)}{2a},\frac{{ (2\phi+1)}
b^2}{3sB}\biggr),\qquad
\lambda^{\mathrm{scad}} \gg\sqrt{\frac{\log(K)}{N}}.
\]

\begin{longlist}[(2)]
\item[(1)] If $\tau_0 \rightarrow0$, then with probability tending
to one,
the LLA--CMA algorithm converges after two LLA iterations and the
LLA--CMA solution (or its one-step version) is equal to the oracle estimator.

\item[(2)] Consider using the LASSO estimator as $\widetilde\boldbeta{}^{(0)}$.
Assume the ir-represent\-able condition in Theorem \ref{thm2}, and pick
the LASSO penalty parameter $\lambda^{\mathrm{lasso}}$ satisfying
\begin{eqnarray*}
\frac{1}{\sqrt{N}} &\ll&\lambda^{\mathrm{lasso}} \ll\min\biggl(\frac{1}{\sqrt{s}}
\rho(s,N),\frac{1}{s}\biggr),
\\
\lambda^{\mathrm{lasso}} &<& \frac{\lambda^{\mathrm{scad}}}{\sqrt{s}}\frac{b}{16}.
\end{eqnarray*}
Then $\tau_0 \rightarrow0$, and the conclusion in (1) holds.
\end{longlist}
\end{cor}
%
%
\begin{remark}\label{remark5}
Part (1) of Corollary \ref{thm5} basically says that any estimator
that converges to $\boldbeta^*$ in probability at a rate faster than
$\lambda^{\mathrm{scad}}$
can be used as the starting value in the LLA--CMA algorithm to find the
oracle estimator with high probability.
Note that such a condition is not very restrictive. Part~(2) of
Corollary \ref{thm5} shows that the LASSO estimator satisfies that condition.
We could also consider using other estimators as the starting value in
the LLA--CMA algorithm. For example, we can use the neighborhood selection
estimator as $\widetilde\boldbeta{}^{(0)}$. Following \citet{RavWaiLaf10}
we assume an ir-representable condition for each of the $K$
neighborhood LASSO-penalized logistic
regression and some other regularity conditions. Then it is not hard to
show that the neighborhood selection
estimator is also a qualified starting value. In this work, we would
like to
faithfully follow the composite likelihood idea and hence prefer to use
the LASSO-penalized composite likelihood estimator as the starting
value in the LLA--CMA algorithm.
\end{remark}

\section{Simulation}\label{sec4}
In this section we use simulation to study
the finite sample performance of the SCAD-penalized
composite likelihood estimator. For comparison, we also include
other two methods: neighborhood selection by LASSO-penalized
logistic regression [\citet{RavWaiLaf10}] and the
LASSO-penalized composite likelihood estimator.

For each coupling coefficient $\beta_{jk}$, the LASSO-penalized
logistic method provides two estimates: $\widehat\beta_{j \mapsto k}$
based on the model for the $j$th dipole and~$\widehat\beta_{k \mapsto
j}$ based on the model for the $k$th dipole. Then we carry\vspace*{2pt}
out two types of neighborhood selections: (i) aggregation by
intersection (NSAI) based
on~$\widehat\beta{}^{\mathrm{NSAI}}_{jk}$,\vadjust{\goodbreak}
and (ii) aggregation by union (NSAU) based on
$\widehat\beta{}^{\mathrm{ NSAU}}_{jk}$, where
\[
\widehat\beta{}^{\mathrm{NSAI}}_{jk}= %
\cases{ 0, &\quad if $\widehat
\beta_{j \mapsto k} \widehat\beta_{k \mapsto j}=0$,
\vspace*{2pt}\cr
\dfrac{\widehat\beta_{j \mapsto k}+\widehat\beta_{k \mapsto j}}{2}, &\quad
otherwise,} %
\]
and
\[
\widehat\beta{}^{\mathrm{NSAU}}_{jk}= %
\cases{ 0, &\quad if $\widehat
\beta_{j \mapsto k}=0 \mbox{ and } \widehat\beta_{k \mapsto j}=0$,
\vspace*{2pt}\cr
\widehat
\beta_{j \mapsto k}, &\quad if $\widehat\beta_{j \mapsto k} \neq0 \mbox{
and }
\widehat
\beta_{k \mapsto j}=0$,
\vspace*{2pt}\cr
\widehat\beta_{k \mapsto j}, &\quad if $\widehat
\beta_{j \mapsto k}=0 \mbox{ and } \widehat\beta_{k \mapsto j} \neq0$,
\vspace*{2pt}\cr
\dfrac{\widehat\beta_{j \mapsto k}+\widehat\beta_{k \mapsto j}}{2}, &\quad if
$\widehat
\beta_{j \mapsto k} \widehat\beta_{k \mapsto j}\neq
0$.}
\]

As suggested by a referee, the relaxed LASSO [\citet{Mei07}] was
used in neighborhood selection to try to improve its estimation accuracy.
In each neighborhood logistic regression model, we first found a subset
model by using the LASSO-penalized
logistic regression. We re-estimated the nonzero coefficients via the
unpenalized logistic regression on the subset model.

BIC has been shown to perform very well for selecting the tuning
parameter of the penalized likelihood estimator
[\citet{WanLiTsa07},
\citet{StaBuhvan10}, \citet{SchBuhvan11}]. We used BIC to
tune all competitors.

%
%
\begin{figure}

\includegraphics{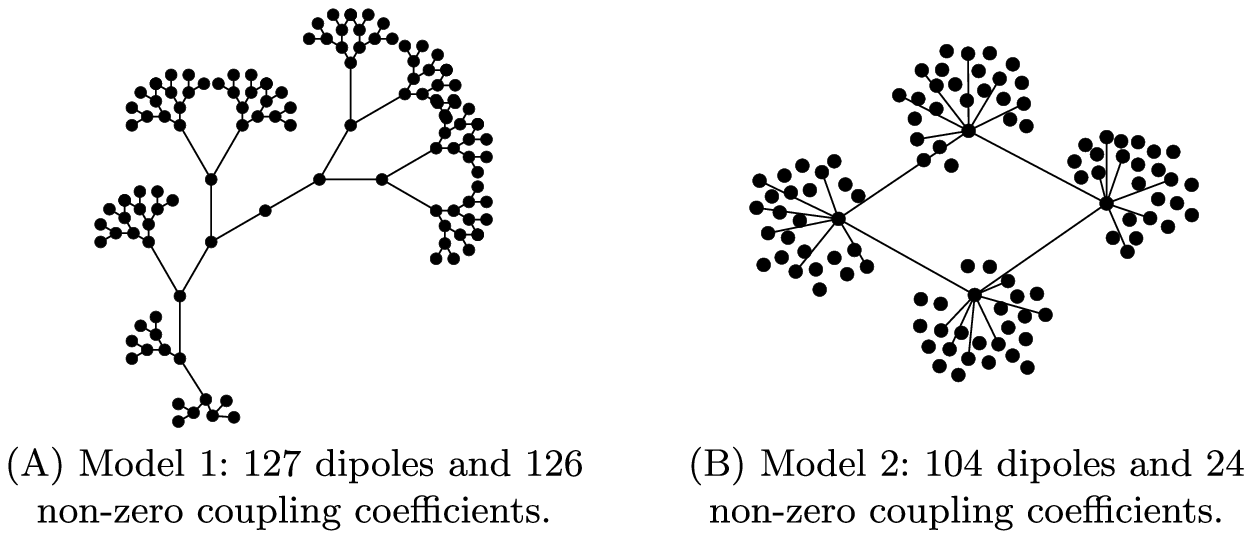}

\caption{Plots of two simulated Ising models.}\label{figur1}
\end{figure}

Two sparse Ising models were considered in our simulation. Their
graphical structure is displayed in Figure \ref{figur1} where solid dots
represent the dipoles, and two dipoles are connected if and only if
their coupling coefficient is nonzero. We generated the nonzero
coupling coefficients as follows. If dipoles~$i$ and $j$ are
connected, we let $\beta_{ij}$ be $t_{ij}s_{ij}$ where $t_{ij}$ is
a random variable following the uniform distribution on $[1,2]$
and $s_{ij}$ is a Bernoulli variable with
$\pr(s_{ij}=1)=\pr(s_{ij}=-1)=0.5$. For each model, we used Gibbs
sampling to generate 100 independent datasets consisting 300
observations. For comparison, we use three measurements: the total
number of discovered edges (NDE), the false discovery rate (FDR)
and mean square errors (MSE).

%
%
\begin{table}
\caption{Comparing different estimators using simulation models 1 and
2 with standard errors in the bracket.
NSAI-relax and NSAU-relax
mean that we use the relaxed LASSO to re-estimate the nonzero
coefficients chosen by neighborhood selection method}\label{tabl1}
\begin{tabular*}{\tablewidth}{@{\extracolsep{\fill}}ld{2.3}d{3.2}cd{2.3}d{3.2}c@{}}
\hline
& \multicolumn{3}{c}{\textbf{Model 1}} &\multicolumn{3}{c@{}}{\textbf{Model 2}} \\[-4pt]
& \multicolumn{3}{c}{\hrulefill} &\multicolumn{3}{c@{}}{\hspace*{1pt}\hrulefill} \\
& \multicolumn{1}{c}{\hspace*{-2pt}\textbf{MSE}}
& \multicolumn{1}{c}{\hspace*{-2pt}\textbf{NDE}}
& \multicolumn{1}{c}{\textbf{FDR}}
& \multicolumn{1}{c}{\hspace*{-2pt}\textbf{MSE}}
& \multicolumn{1}{c}{\hspace*{-1pt}\textbf{NDE}}
& \multicolumn{1}{c@{}}{\textbf{FDR}} \\
\hline
NSAI
&22.96 &138.9 &0.09 &8.16 &26.8 &0.16 \\
&(0.18) &(0.4) &(0.01) &(0.12) &(0.2) &(0.01) \\
NSAU
&17.34 &197.3 &0.36 &6.38 &39.7 &0.39 \\
&(0.14) &(1.0) &(0.01) &(0.16) &(0.5) &(0.01) \\
LASSO
&21.33 &332.5 &0.62 &12.19 &117.1 &0.79 \\
&(0.13) &(3.8) &(0.04) &(0.12) &(3.0) &(0.05) \\
SCAD1
&2.86 &145.0 &0.12 &5.64 &30.0 &0.22 \\
&(0.10) &(2.4) &(0.01) &(0.17) &(1.8) &(0.02) \\
SCAD2
&2.43 &129.2 &0.07 &4.41 &26.1 &0.17 \\
&(0.05) &(0.5) &(0.01) &(0.13) &(0.7) &(0.02) \\
SCAD2$^{**}$
&2.42 &128.6 &0.06 &4.39 &25.7 &0.16 \\
&(0.05) &(0.5) &(0.01) &(0.13) &(0.6) &(0.02) \\
NSAI-relax
&8.23 &138.9 &0.09 &6.34 &26.8 &0.16 \\
&(0.13) &(0.4) &(0.01) &(0.09) &(0.2) &(0.01) \\
NSAU-relax
&4.44 &197.3 &0.36 &5.67 &39.7 &0.39 \\
&(0.10) &(0.4) &(0.01) &(0.10) &(0.5) &(0.01) \\
\hline
\end{tabular*}
\end{table}

Based on Table \ref{tabl1}, we make the following interesting observations:

\begin{itemize}
\item NSAU, while selecting larger models than NSAI, provides more
accurate estimation. Neighborhood selection outperforms
the LASSO-penalized composite likelihood estimator.

\item Note that SCAD2$^{**}$ has the smallest MSE in both models.
SCAD2$^{**}$ and SCAD2 gave almost identical results, and their
improvement over SCAD1 is statistically significant.
All three SCAD solutions perform much better than the LASSO for fitting
penalized composite likelihood in terms of estimation
and selection.

\item The SCAD solutions and NSAI have similar model selection
performance, but the
SCAD is substantial better in estimation. Using the relaxed LASSO can
improve the estimation accuracy of neighborhood selection methods,
but their improved MSEs are still significantly higher than those of
SCAD2 and SCAD2$^{**}$.

\end{itemize}

In Table \ref{tabl2} we compare the run times of the three methods. LASSO-CGA
denotes the coordinate gradient ascent algorithm for computing the
LASSO estimator. The computing time is about five times longer than
that used by the CMA algorithm. Compared to the LASSO case, the run
time for fitting the SCAD model is doubled or tripled, but it is still
very manageable for the high-dimensional data.

%
%
\begin{table}
\caption{Total time (in seconds) for computing solutions at 100
penalization parameters, averaged over 3 replications. Timing was
carried out on a laptop with an Intel Core 1.60 GHz processor. LASSO-CGA
denotes a coordinate gradient ascent algorithm for computing the
LASSO-penalized composite likelihood estimator. The timing of SCAD1,
SCAD2 and SCAD2$^{**}$ includes the timing for~computing the starting value}
\label{tabl2}
\begin{tabular*}{\tablewidth}{@{\extracolsep{\fill}}lcccccc@{}}
\hline
& \textbf{Neighborhood} & & &\\
$\bolds{(N,p)}$ & \textbf{selection} & \textbf{LASSO}
& \textbf{SCAD1} & \textbf{SCAD2} & \textbf{SCAD2}$\bolds{^{**}}$ &
\textbf{LASSO-CGA} \\
\hline
Model 1& 51.1 & 32.7 &67.9 & 84.7 & 95.1 & 179.8 \\
$(300,7875)$ & & && & \\[4pt]
Model 2& 29.8 & 16.0 &34.8 & 42.6 & 51.2 & \hphantom{0}89.6\\
$(300,5356)$ & & & & & \\
\hline
\end{tabular*}
\end{table}

\section{Stanford HIV drug resistance data}\label{sec5}

We also illustrate our methods in a real example using a HIV
antiretroviral therapy (ART) susceptibility dataset obtained from
the Stanford HIV drug resistance database. Details of the database
and related data sets can be found in \citet{Rheetal06}. The data
for analysis consists of virus mutation information at 99 protease
residues (sites) for $N = 702$ isolates from the plasma of
HIV-1-infected patients. This dataset has been previously used in
\citet{Rheetal06} and \citet{WuCaiLin10} to study the
association between protease mutations and susceptibility to ART
drugs.

A well recognized problem with current ART treatment such as PIs
for treating HIV is that individuals who initially respond to
therapy may develop resistance to it due to viral mutations. HIV-1
protease plays a key role in the late stage of viral replication
and its ability to rapidly acquire a variety of mutations in
response to various PIs confers the enzyme with high resistance to
ARTs. A high cooperativity has been observed among drug-resistant
mutations in HIV-1 protease [\citet{OhtSchFre03}]. The sequence
data retrieved from treated patients is likely to include
mutations that reflect cooperative effects originating from late
functional constraints, rather than stochastic evolutionary noise
[\citet{Atcetal00}]. However, the molecular mechanisms of drug
resistance is yet to be elucidated. It is thus of great interest
to study inter-residue couplings which might be relevant to
protein structure or function and thus could potentially shed
light on the mechanisms of drug resistance. We apply the proposed
method to the protease sequence data to investigate such
inter-residue contacts. Our analysis only included $K=79$ of the
99 residues that contain mutations.

We split the data into a training set with 500 data and a test set
with 202 data. Model fitting and selection were done on the
training set\vspace*{1pt} and the test data were used to compare the model
errors. For a given estimate $\widehat\boldbeta$ obtained from
the training set, its model error
is gauged by the value of composite likelihood evaluated on the
test set, that is,
\[
\operatorname{ME}(\widehat\boldbeta)=-\ell^{\mathrm{test}}_c(\widehat
\boldbeta) =- \frac{1}{202}\sum^{202}_{n=1}
\sum^{79}_{j=1}\log\bigl(
\theta_{jn}(\widehat\boldbeta)\bigr).
\]

%
%
\begin{figure}

\includegraphics{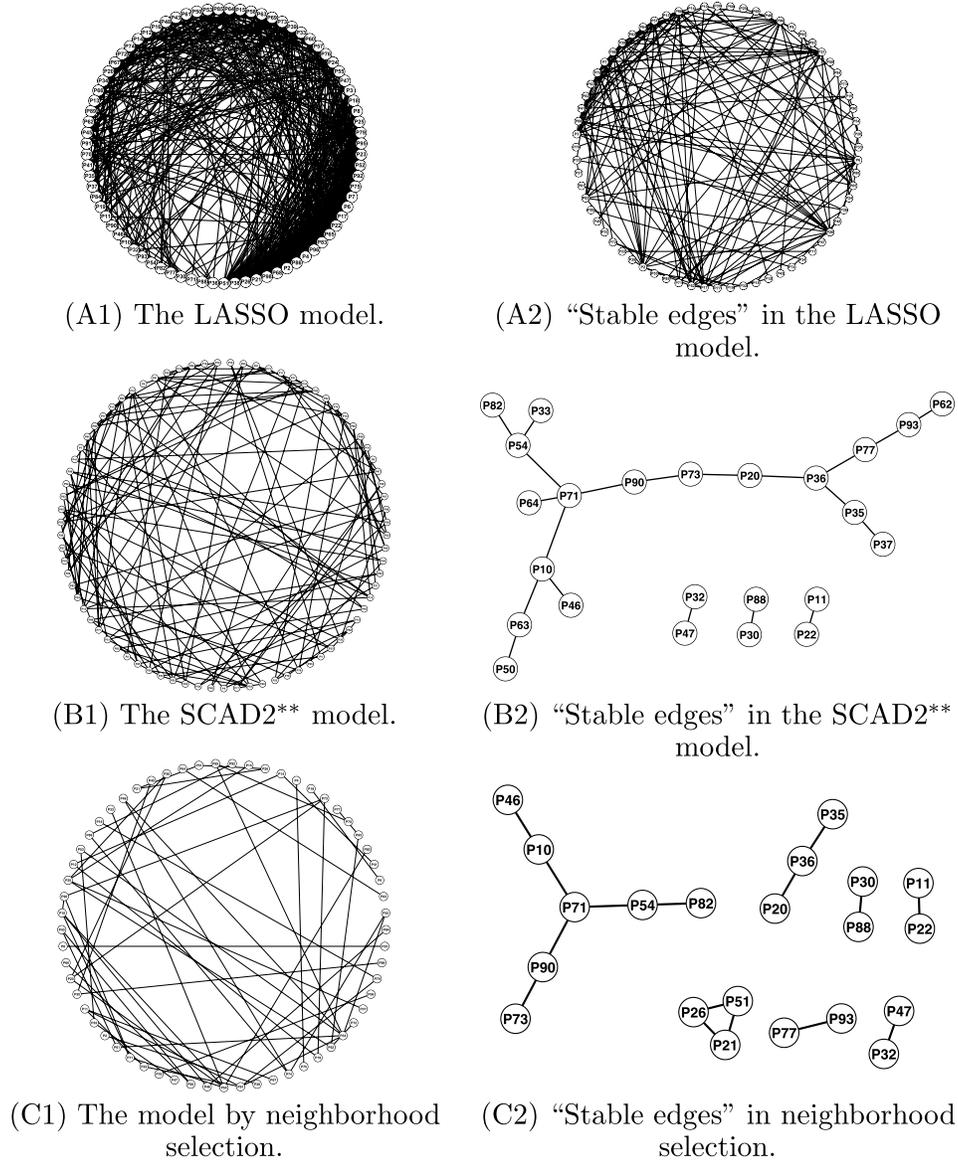}

\caption{Shown in the left three panels \textup{(A1)}, \textup{(B1)},
\textup{(C1)} are the
selected models by BIC. The right three panels \textup{(A2)},
\textup{(B2)}, \textup{(C2)} show the
stability selection results using $\pi_{\mathrm{thr}}=0.9$.}\label{figur2}
\end{figure}

We report the analysis results in Table \ref{tabl3}. There are total
$3081$ coupling coefficients to be estimated. Graphical presentations
of the selected models are shown in Figure \ref{figur2}. Note that SCAD2 and
SCAD2$^{**}$ again gave almost identical results and performed better
SCAD1. We also performed stability selection [\citet{MeiBuh10}] on
each method to find ``stable edges.'' A remarkable property of
%
%
\begin{table}
\caption{Application to HIVRT data. NSE is the number of ``stable
edges.'' $E[V]$ is the expected number of falsely selected edges. Its
upper bounds were computed by Theorem 1 in Meinshausen~and~B{\"u}hlmann
(\protect\citeyear{MeiBuh10})}
\label{tabl3}
\begin{tabular*}{\tablewidth}{@{\extracolsep{\fill}}ld{2.2}d{3.2}d{3.2}d{3.2}d{3.2}d{3.2}@{}}
\hline
& \multicolumn{1}{c}{\textbf{NSAI}} & \multicolumn{1}{c}{\textbf{NSAU}}
& \multicolumn{1}{c}{\textbf{LASSO}} & \multicolumn{1}{c}{\textbf{SCAD1}}
& \multicolumn{1}{c}{\textbf{SCAD2}} & \multicolumn{1}{c@{}}{\textbf{SCAD2}$\bolds{^{**}}$}\\
\hline
NDE & 57 & 305 & 631 & 101 & 141 & 132 \\
ME & 26.38 & 36.34 & 18.35 & 18.30 & 16.76 & 16.74 \\
[4pt]
\multicolumn{7}{@{}c@{}}{Stability selection} 
\\
[4pt]
NSE ($\pi_{\mathrm{thr}}=0.9$)& 15 & 63 & 160 & 17 & 20 & 20 \\
$E[V]$ & \multicolumn{1}{c}{$\le3.2$} & \multicolumn{1}{c}{$\le48$}
& \multicolumn{1}{c}{$\le$147.5} & \multicolumn{1}{c}{$\le$4.3}
& \multicolumn{1}{c}{$\le$8.0} &
\multicolumn{1}{c@{}}{$\le$7.2} \\
\hline
\end{tabular*}
\end{table}
stability selection is that under some suitable conditions stability
selection achieves finite sample control over the expected number of
false discoveries in the set of ``stable edges.'' We use the SCAD
selector to explain the stability selection procedure. We took a~random
subsample of size 250 and fitted the SCAD model. The process was
repeated 100 times. On average, SCAD1 selected $103.1$ edges, SCAD2
selected $140.7$ edges and SCAD2$^{**}$ chose $133.4$ edges. For each
coefficient~$\beta_{jk}$ we computed its frequency of being selected,
denoted by $\widehat\Pi_{jk}$. The set of ``stable edges'' is defined as
$\{(k,j)\dvtx\widehat{\Pi}_{kj}>\pi_{\mathrm{thr}}\}$. In Table \ref{tabl3}, we
report the results using the threshold $\pi_{\mathrm{thr}}=0.9$, as suggested by
\citet{MeiBuh10}. Stability selection found 17 edges
in the SCAD1. SCAD2 and SCAD2$^{**}$ selected the same $20$ stable
edges. By Theorem 1 in \citet{MeiBuh10}, among these
17 stable edges selected by SCAD1, the expected number of false
discoveries is no greater than $4.3$, and among the 20 stable edges
selected by SCAD2 or SCAD2$^{**}$, the expected number of false
discoveries is at most $7.2$. Likewise, we did stability selection with
the LASSO selector and neighborhood selection, and the results are
reported in Table \ref{tabl3} as well. Figure \ref{figur2} shows the ``stable
edges'' by stability selection. We see that the computed upper bounds
are very useful for the SCAD selector and NSAI and not so informative
for the LASSO selector and NSAU. Interestingly, both NSAI and SCAD
suggest there are about $12$ true discoveries by stability selection.
In fact, we found that NSAI and SCAD1 have 11 ``stable edges'' in
common, and NSAI and SCAD2 (or SCAD2$^{**}$) have 12 ``stable edges''
in common.

These results are consistent with some of the previous findings.
For example, it has long been known that co-substitutions at
residues 30 and 88 are most effective in reducing the
susceptibility of nelfinavir [\citet{LiuEyaBah08}]. Among the
top 30
most common drug resistance mutations [\citet{Rheetal04}], 7~of
those had a joint mutation at residues 54 and 82, the joint
mutation at residues 88 and 30 was the second most common
mutation among all. A co-mutation at residues 54, 82 and 90 was
associated with high resistance to multiple drugs and an
additional co-mutation at 46 was associated with an even higher
level of resistance. It is interesting to note that using a larger
set of isolates from treated HIV patients, \citet{Wuetal03}
reported (54, 82), (32, 47), (73, 90) as the three most highly
correlated pairs. All these three pairs showed up as the stable
edges in our analysis. Mutation at residue 71, often described
as a compensatory or accessory mutation, has been reported as a
critical mutation which appears to improve virus growth and
contribute to resistance phenotype [\citet{Maretal},
\citet{Tisetal95}, \citet{MuzRosFre03}]. Accessory mutations
contribute to resistance only when present with a mutation in the
substrate cleft or flap or at residue 90 [\citet{Wuetal03}]. The
stable edges connect this accessory mutation with residues 90 and
54 (a flap residue), as well as with another flap residue at 46
through residue~10.

\begin{appendix}\label{app}
\section*{Appendix: Technical proofs}
Before presenting the proof, we first define some useful
quantities. The score functions of the negative composite likelihood
($-\ell^{(j)}$) and the Hessian
matrices are defined as follows:
\begin{eqnarray*}
\psi^{(j)}_k&=&-\frac{\partial\ell^{(j)}(\boldbeta^{(j)})}{\partial
\beta_{jk}}=\frac{1}{N}\sum
^N_{n=1}x_{jn}x_{kn}(
\theta_{jn}-1),\qquad k \neq j,
\\
H^{(j)}_{k_1,k_2}&=&-\frac{\partial^2
\ell^{(j)}(\boldbeta^{(j)})}{\partial\beta_{jk_1}\,\partial
{\beta_{jk_2}}}=\frac{1}{N}\sum
^N_{n=1}x_{k_1 n}x_{k_2
n}(1-
\theta_{jn})\theta_{jn},\qquad k_1,k_2 \neq
j.
\end{eqnarray*}
Similarly, let $\psi$ be the score function of $-\ell_c$ such that
$\psi_{(jk)}=\frac{\partial-\ell_c(\boldbeta)}{\partial
\beta_{jk}}$ for $1 \le j<k \le K$. By definition we have the
following identities: $ \psi_{(jk)}=\psi^{(j)}_k+\psi^{(k)}_j $. In
what follows we write $\psi^*=\psi(\boldbeta^*)$.

\begin{pf*}{Proof of Theorem \ref{thm1}}
We first prove part (1).

Consider $ V(\boldalpha_{\mathcal{A}})=
-\ell_c(\boldbeta_{\mathcal{A}}^* +
d_N\boldalpha_{\mathcal{A}})+\ell_c(\boldbeta_{\mathcal{A}}^*) $
and its minimizer is $
\widetilde{\boldalpha}{}^{\mathrm{hmle}}_{\mathcal{A}}=\frac
{1}{d_N}(\widetilde{\boldbeta}{}^{\mathrm{hmle}}_{\mathcal{A}}
- \boldbeta_{\mathcal{A}}^*) $. By definition, $
V(\widetilde{\boldalpha}{}^{\mathrm{hmle}}_{\mathcal{A}})\le V(\boldzero)=0
$. Fix any $R>0$ and consider any $\boldalpha_{\mathcal{A}}$
satisfying $ \Vert\boldalpha_{\mathcal{A}} \Vert_2=R $. Using
Taylor's expansion, we know that, for some $t\in[0,1]$ and $
\boldbeta(t)=\boldbeta^*_{\mathcal{A}}+td_N\boldalpha_{\mathcal{A}}$,
%
%
\begin{eqnarray}
\label{expan}\quad
V(\boldalpha_{\mathcal{A}}) &=& d_N
\boldalpha^T_{\mathcal{A}}\psi^*_{\mathcal{A}}+
\tfrac{1}{2}d^2_N \boldalpha^T_{\mathcal{A}}
H^*_{\mathcal{AA}} \boldalpha_{\mathcal{A}} \nonumber\\
&&{}+ \tfrac{1}{2}d^2_N
\boldalpha^T_{\mathcal{A}} \bigl[H_{\mathcal{AA}}\bigl(\boldbeta(t)
\bigr) -H^*_{\mathcal{AA}} \bigr] \boldalpha_{\mathcal{A}}
\\
&\equiv&T_1+T_2+T_3.\nonumber
\end{eqnarray}
Note that $E[\psi^*_{\mathcal{A}}]=0$ and
$\Vert\psi^*_{\mathcal{A}}\Vert_{\infty}\le2$. By the Cauchy--Schwarz
inequality, $ \vert\boldalpha^T_{\mathcal{A}}\psi^*_{\mathcal{A}}
\vert\le2\sqrt{s}R$. Using Hoeffding's inequality, we have
%
%
\begin{equation}
\label{T1} \Pr( T_1 \ge-d_N \varepsilon) \le\exp\biggl(-
\frac{N\varepsilon^2}{8sR^2}\biggr).
\end{equation}
For the second term, we first have $ T_2 \ge\frac{d^2_N}{2}
\lambda_{\min}(H^*_{\mathcal{AA}})R^2. $ Each entry of $H^*$ is
between $-\frac{1}{2}$ and $\frac{1}{2}$. Thus Hoeffding's
inequality and the union bound yield
\[
\Pr\biggl(\bigl\Vert H^{(N)}_j-H_j
\bigr\Vert_F^2\ge\frac{b^2}{4}\biggr) \le2s^2
\exp\biggl(-N\frac{b^2}{2s^2}\biggr).
\]
So by the inequality $ \lambda_{\min}(H^*_{\mathcal{AA}}) \ge
b-\Vert H^*_{\mathcal{AA}}-E[H^*_{\mathcal{AA}}]\Vert_{F}$, we
have
%
%
\begin{equation}
\label{T2} \Pr\bigl(T_2 \ge d^2_NbR^2/4
\bigr) \ge1-2s^2\exp\biggl(-\frac{Nb^2}{2s^2}\biggr).
\end{equation}
For $|T_3|$, let $
\lambda_{\max}(\frac{1}{N}\sum^N_{n=1}\boldx_{{\mathcal A} n}
\boldx_{{\mathcal A} n}^T)=B_N$. Define $\bar
\eta_{jn}(\boldbeta)=\theta_{jn} (1-\theta_{jn} )
(2\theta_{jn}-1 )$. Using the mean value theorem, we have
that, for some $t^{\prime}\in[0,t]$ and $
\boldbeta(t')=\boldbeta^*_{\mathcal{A}}+t'd_N\boldalpha_{\mathcal{A}}$,
%
%
\begin{eqnarray}
\label{bT3}
|T_3| &=& \frac{d^3_N}{2}\Biggl|{\frac{1}{N}\sum
_n}\sum_{j=1}^K
\mathop{\sum_{k_1\neq j}}_{
k_2\neq j}
\alpha_{jk_1}\alpha_{jk_2}x_{k_1n}x_{k_2n}
t'\bar\eta_{jn}\bigl(\boldbeta\bigl(t'\bigr)
\bigr) \biggl(\sum_{{ k' \neq j} }\alpha_{jk'}x_{jn}x_{k'n}
\biggr)\Biggr|
\nonumber\hspace*{-35pt}\\[-8pt]\\[-8pt]
& \le&
\frac{d^3_N}{2} \biggl(\frac{\sqrt{sR^2}}{4}\biggr) \cdot\biggl(2 B_N
\sum_{(j,k) \in\mathcal{A}} \alpha^2_{jk}
\biggr)=\frac{d^3_NB_N}{4}\sqrt{s}R^3.
\nonumber\hspace*{-35pt}
\end{eqnarray}
In the last step we have used $|\bar\eta_{jn}(\boldbeta(t'))|\le
\frac{1}4$ for any $j$ and $\boldalpha_{\mathcal{A}^c}=0$. Moreover,
$ B_N \le B +\Vert
\frac{1}{N}\sum^N_{n=1}\boldx_{{\mathcal A} n} \boldx_{{\mathcal
A} n}^T -E[\boldx_{{\mathcal A}} \boldx_{{\mathcal A}}^T]\Vert_F.
$ Since\vspace*{1pt} $x_{jn}=\pm1$, we apply Hoeffding's inequality and the
union bound to obtain the following probability bound:
\[
\Pr\Biggl(\Biggl\Vert\frac{1}{N}\sum^N_{n=1}
\boldx_{{\mathcal A} n} \boldx_{{\mathcal A} n}^T -E\bigl[
\boldx_{{\mathcal A}} \boldx_{{\mathcal A}}^T\bigr]
\Biggr\Vert_F \ge B/2\Biggr) \le2s^2\exp\biggl(-
\frac{NB^2}{8s^2}\biggr),
\]
which leads to
%
%
\begin{equation}
\label{T3} \Pr\biggl(|T_3| \le\frac{3d^3_NB}{8}\sqrt{s}R^3
\biggr) \ge1-2s^2\exp\biggl(-\frac{NB^2}{8s^2}\biggr).
\end{equation}
Taking $R<\frac{b}{3B}\frac{\sqrt{N}}{s}$ and combining (\ref{T1})
(\ref{T2}) and (\ref{T3}), we have
\[
T_1+T_2+T_3 \ge\frac{bR^2}{8}d^2_N-
\frac{3B}{8}R^3d^3_N\sqrt{s}>0
\]
with probability at least $1-\tau_1$. Thus, the convexity of $V$
implies that
\[
\Pr\biggl(\bigl\Vert\widetilde\boldbeta{}^{\mathrm{hmle}}_{\mathcal
A}-
\boldbeta^*_{\mathcal A} \bigr\Vert_2 \leq\sqrt{\frac{s}{N}} R
\biggr) \ge1-\tau_1.
\]

We now prove part (2). First, we show that if $\min_{(j,k) \in
\mathcal{A}} |\widetilde\beta{}^{\mathrm{hmle}}_{jk}|>a\lambda$ and $\Vert{
\psi_{{\mathcal A}^c}(\widehat\boldbeta{}^{\mathrm{oracle}})}\Vert_{\infty}
\leq\lambda$, then $\widehat{\boldbeta}{}^{\mathrm{oracle}}$ is a local
maximizer of $\ell_c(\boldbeta)-
\sum_{(j,k)}P_{\lambda}(|\beta_{jk}|)$. To see that, consider
a small ball of radius $t$ with $\widehat{\boldbeta}{}^{\mathrm{oracle}}$
being the center. Let $\boldbeta$ be any point in the ball. So
$\Vert\boldbeta-\widehat\boldbeta{}^{\mathrm{oracle}}\Vert_2 \leq t$.
Clearly, for a sufficiently small $t$ we have $\min_{(j,k) \in
\mathcal{A}} | \beta_{jk}|>a\lambda$ and $\max_{(j,k) \in
\mathcal{A}^c} | \beta_{jk}|<\lambda$. By Taylor's expansion we
have
\begin{eqnarray*}
&&\biggl\{-\ell_c(\boldbeta)+\sum_{(j,k)}P_{\lambda}\bigl(|
\beta_{jk}|\bigr)\biggr\}-\biggl\{ -\ell_c\bigl(\widehat
\boldbeta{}^{\mathrm{oracle}}\bigr)+\sum_{(j,k)}P_{\lambda}
\bigl(\bigl|\widehat\beta{}^{\mathrm{oracle}}_{jk}\bigr|\bigr)\biggr\}
\\
&&\qquad= \bigl(\boldbeta_{\mathcal A}-\widetilde\boldbeta{}^{\mathrm{hmle}}
\bigr)^T { \psi_{{\mathcal A}^c}\bigl(\widehat\boldbeta{}^{\mathrm{oracle}}
\bigr)}+\frac
{1}{2}\bigl(\boldbeta-\widehat\boldbeta{}^{\mathrm{oracle}}
\bigr)^TH\bigl(\boldbeta'\bigr) \bigl(\boldbeta-\widehat
\boldbeta{}^{\mathrm{oracle}}\bigr)
\\
&&\qquad\quad{}+\sum_{(j,k) \in\mathcal{A}^c}\lambda|\beta_{jk}|
\\
&&\qquad\ge \sum_{(j,k) \in\mathcal{A}^c}\bigl(\lambda-\bigl|\psi_{(jk)}
\bigl(\widehat\boldbeta{}^{\mathrm{oracle}}\bigr)\bigr|\bigr) |\beta_{jk}| \ge0.
\end{eqnarray*}

A probability bound for the event of $\min_{(j,k) \in\mathcal
{A}} |\widetilde\beta{}^{\mathrm{hmle}}_{jk}|>a\lambda$ is given by
%
%
\begin{eqnarray}
\label{minimum} &&\Pr\Bigl(\min_{(j,k) \in\mathcal{A}} \bigl|\widetilde\beta
^{\mathrm{hmle}}_{jk}\bigr|>a
\lambda\Bigr)
\nonumber
\\
&&\qquad \ge\Pr\Biggl(\bigl\Vert\widetilde\boldbeta{}^{\mathrm{hmle}}_{\mathcal
A}-
\boldbeta^*_{\mathcal A} \bigr\Vert_2 \leq\sqrt{\frac{s}{N}}R_*
\Biggr)
\\
&&\qquad\ge 1-\exp\biggl(-R^{2}_*\frac{b^2}{8^3}\biggr)-2s^2
\exp\biggl(-\frac{N}{s^2}\frac
{b^2}{2}\biggr)-2s^2\exp
\biggl(-\frac{N}{s^2}\frac{B^2}{8}\biggr).\nonumber
\end{eqnarray}
Now consider $ \Pr(\Vert{\psi_{{\mathcal A}^c}(\widehat\boldbeta
{}^{\mathrm{oracle}})}\Vert_{\infty} < \lambda). $
There exists some $t \in[0,1]$ such that
%
%
\begin{equation}
\label{scadp2eq1} \psi\bigl(\widehat\boldbeta{}^{\mathrm{oracle}}\bigr)=\psi\bigl(
\boldbeta^*\bigr)+H^{*} \bigl(\widehat\boldbeta{}^{\mathrm{oracle}}-
\boldbeta^*\bigr)+r,
\end{equation}
where $r=(H(\boldbeta^*+t(\widehat
\boldbeta{}^{\mathrm{oracle}}-\boldbeta^*))-H^{*})(\widehat
\boldbeta{}^{\mathrm{oracle}}-\boldbeta^*)$. Note $\psi_{{\mathcal A}}(\widehat
\boldbeta{}^{\mathrm{oracle}})=0$, so
\[
\widetilde\boldbeta_{{\mathcal A}}-\boldbeta^*_{{\mathcal
A}}=
\bigl(H^*_{\mathcal AA}\bigr)^{-1}(-\psi_{\mathcal A}-r_{\mathcal A}).
\]
Then $\Vert{\psi_{{\mathcal A}^c}(\widehat\boldbeta{}^{\mathrm{oracle}})}\Vert
_{\infty} \leq\lambda$ becomes
\[
\bigl\Vert H^{*}_{\mathcal A^cA}\bigl(H^{*}_{\mathcal
AA}
\bigr)^{-1}(-\psi_{\mathcal A}-r_{\mathcal A})+\psi_{\mathcal
A^c}+r_{\mathcal A^c}
\bigr\Vert_{\infty}\le\lambda,
\]
which is guaranteed if
\[
{ \bigl(\bigl\Vert H^{*}_{\mathcal
A^cA}\bigl(H^{*}_{\mathcal AA}
\bigr)^{-1} \bigr\Vert_{\infty}+1\bigr) }\bigl(\Vert\psi
\Vert_{\infty}+\Vert r \Vert_{\infty}\bigr) \le\lambda.
\]
Therefore we
have a simple lower bound for $ \Pr(\Vert{\psi_{{\mathcal
A}^c}(\widehat\boldbeta{}^{\mathrm{oracle}})}\Vert_{\infty} \le\lambda)$.
\begin{eqnarray*}
&&\Pr\bigl(\bigl\Vert{ \psi_{{\mathcal A}^c}\bigl(\widehat
\boldbeta{}^{\mathrm{oracle}}\bigr)}\bigr\Vert_{\infty} \le\lambda\bigr)
\\
&&\qquad>1-\Pr\bigl(\bigl\Vert H^{*}_{\mathcal A^cA}\bigl(H^{*}_{\mathcal
AA}
\bigr)^{-1}\bigr\Vert_{\infty}>2\phi\bigr)
-\Pr\biggl(\Vert\psi
\Vert_{\infty}>\frac{\lambda}{{ 4\phi+2}}\biggr)\\
&&\qquad\quad{}-\Pr\biggl(\Vert r
\Vert_{\infty}>\frac{\lambda}{{ 4\phi+2}}\biggr).
\end{eqnarray*}
Using Hoeffding's inequality and the union bound, we have
%
%
\begin{equation}
\label{scadbdelta5} \Pr\biggl(\Vert\psi\Vert_{\infty} \le\frac{\lambda
}{{ 4\phi+2}}
\biggr) \ge1-K^2\exp\biggl(-\frac{N\lambda^2}{{ 128(\phi+
{1}/2)^2}}\biggr).
\end{equation}
Write $\boldalpha=\widetilde\boldbeta{}^{\mathrm{hmle}}-\boldbeta^*$, and thus
$\boldalpha_{\mathcal A^c}=0$. By the
mean value theorem, we have a bound for $r_{(jk)}$:
\begin{eqnarray*}
|r_{(jk)}| &=& \Biggl|\frac{1}{N}\sum^N_{n=1}
\sum_{k_2\neq j} \sum_{{k'\neq j}}
x_{kn}x_{jn}x_{k_2n}x_{k'n}
\alpha_{jk_2}\alpha_{jk'} t'\bar
\eta_{jn}\bigl(\boldbeta\bigl(t'\bigr)\bigr)
\nonumber
\\
&&\hspace*{2pt}{}+\frac{1}{N}\sum^N_{n=1} \sum
_{j_2\neq k} \sum_{{j'\neq k}}
x_{jn}x_{kn}x_{j_2n}x_{j'n}
\alpha_{kj_2}\alpha_{kj'} t'\bar
\eta_{kn}\bigl(\boldbeta\bigl(t'\bigr)\bigr)\Biggr|
\nonumber
\\
& \le& B_N\cdot\bigl\Vert\widetilde\boldbeta_{\mathcal A}-
\boldbeta^*_{\mathcal A} \bigr\Vert^2_2.
\end{eqnarray*}
In the last step we have used $|\bar\eta_{jn}(\boldbeta(t'))|\le
\frac{1}4$ for any $j$ and $\boldalpha_{\mathcal A^c}=0$. Moreover,
recall that
\[
B_N \le B +\Biggl\Vert\frac{1}{N}\sum^N_{n=1}
\boldx_{{\mathcal A} n} \boldx_{{\mathcal A} n}^T -E\bigl[
\boldx_{{\mathcal A}} \boldx_{{\mathcal A}}^T\bigr]
\Biggr\Vert_F.
\]
Thus
%
%
\begin{eqnarray}
\label{scadbdelta6}
\Pr\biggl(\Vert r \Vert_{\infty} < \frac{\lambda
}{{ 4\phi+2}}
\biggr) &\ge&1-{\exp\biggl(\frac{- N\lambda}{3B(2\phi+1)s}\frac
{b^2}{8^3}\biggr)}-2s^2
\exp\biggl(\frac{-Nb^2}{2s^2}\biggr)\nonumber\hspace*{-35pt}\\[-8pt]\\[-8pt]
&&{}-2s^2\exp\biggl(\frac{-NB^2}{8s^2}
\biggr).\nonumber\hspace*{-35pt}
\end{eqnarray}

For notation convenience define $c=\Vert{(E[H^{*}_{\mathcal
AA}])^{-1}} \Vert_{\infty} \leq\sqrt{s} \Vert{(E[H^{*}_{\mathcal
AA}])^{-1}} \Vert_{2}$ and
\begin{eqnarray*}
\delta&=&\bigl\Vert H^{*}_{\mathcal A^cA}\bigl(H^{*}_{\mathcal AA}
\bigr)^{-1} -E\bigl[H^{*}_{\mathcal A^cA}\bigr]\bigl(E
\bigl[H^{*}_{\mathcal
AA}\bigr]\bigr)^{-1}
\bigr\Vert_{\infty},
\\
\delta_1&=&\bigl\Vert\bigl(H^{*}_{\mathcal AA}
\bigr)^{-1}-{\bigl(E\bigl[H^{*}_{\mathcal
AA}\bigr]
\bigr)^{-1}} \bigr\Vert_{\infty},
\\
\delta_2&=&\bigl\Vert H^{*}_{\mathcal AA}-E
\bigl[H^{*}_{\mathcal AA}\bigr] \bigr\Vert_{\infty},
\\
\delta_3&=&\bigl\Vert H^{*}_{\mathcal
A^cA}-E
\bigl[H^{*}_{\mathcal A^cA}\bigr] \bigr\Vert_{\infty}.
\end{eqnarray*}
Then by definition
\begin{eqnarray*}
\delta&=&\bigl\Vert\bigl(H^{*}_{\mathcal A^cA}-E\bigl[H^{*}_{\mathcal
A^cA}
\bigr]\bigr) \bigl(\bigl(H^{*}_{\mathcal AA}\bigr)^{-1}-{
\bigl(E\bigl[H^{*}_{\mathcal
AA}\bigr]\bigr)^{-1}}\bigr)
\nonumber
\\
&&\hspace*{4.5pt}{}+E\bigl[H^{*}_{\mathcal A^cA}\bigr]{\bigl(E\bigl[H^{*}_{\mathcal
AA}
\bigr]\bigr)^{-1}}\bigl(-H^{*}_{\mathcal AA}+E
\bigl[H^{*}_{\mathcal AA}\bigr]\bigr) \bigl(H^{*}_{\mathcal
AA}
\bigr)^{-1}
\nonumber
\\
&&\hspace*{88pt}{}+\bigl(H^{*}_{\mathcal A^cA}-E\bigl[H^{*}_{\mathcal A^cA}
\bigr]\bigr) {\bigl(E\bigl[H^{*}_{\mathcal AA}\bigr]
\bigr)^{-1}}\bigr\Vert_{\infty}
\nonumber
\\
& \le& \delta_3 \delta_1+\phi\delta_2\bigl\Vert
\bigl(H^{*}_{\mathcal
AA}\bigr)^{-1} \bigr\Vert_{\infty}+
\delta_3 c
\nonumber
\\
& \le& \delta_3 \delta_1+ \phi(c+\delta_1)
\delta_2+ \delta_3 c.
\end{eqnarray*}
Note that
\begin{eqnarray*}
\delta_1&=& \bigl\Vert\bigl(H^{*}_{\mathcal AA}
\bigr)^{-1}\bigl(E\bigl[H^{*}_{\mathcal
AA}
\bigr]-H^{*}_{\mathcal AA}\bigr){\bigl(E\bigl[H^{*}_{\mathcal AA}
\bigr]\bigr)^{-1}} \bigr\Vert_{\infty}
\\
& \le& \bigl\Vert\bigl(H^{*}_{\mathcal AA}\bigr)^{-1}
\bigr\Vert_{\infty} \cdot\bigl\Vert E\bigl[H^{*}_{\mathcal AA}
\bigr]-H^{*}_{\mathcal AA} \bigr\Vert_{\infty} \cdot\bigl\Vert{\bigl(E
\bigl[H^{*}_{\mathcal AA}\bigr]\bigr)^{-1}}
\bigr\Vert_{\infty}
\\
& \le& (\delta_1+c) \delta_2 c.
\end{eqnarray*}
Hence as long as $\delta_2c<1$ we have $ \delta_1 \le
\frac{\delta_2c^2}{1-\delta_2 c} $ and $\delta\le
(\delta_3+\phi\delta_2)\frac{c}{1-\delta_2c}. $
%
%
\begin{eqnarray}
\label{scadbdelta7} \Pr\biggl(\delta_2<\frac{1}{4c}\biggr) &
\ge& 1- \Pr\biggl(\bigl\Vert H^{*}_{\mathcal A^cA}-E\bigl[H^{*}_{\mathcal A^cA}
\bigr]\bigr\Vert_{\max}>\frac{1}{4cs}\biggr)
\nonumber\\[-8pt]\\[-8pt]
&\ge&1-2s^2\exp\biggl(-\frac{N}{8c^2s^2}\biggr),
\nonumber\\
%
%
\label{scadbdelta8} \Pr\biggl(\delta_3<\frac{\phi}{2c}\biggr) &
\ge& 1- \Pr\biggl(\bigl\Vert H^{*}_{\mathcal A^cA}-E\bigl[H^{*}_{\mathcal A^cA}
\bigr]\bigr\Vert_{\max}>\frac{\phi}{4cs}\biggr)
\nonumber\\[-8pt]\\[-8pt]
&\ge&1-K^2s\exp\biggl(-\frac{N\phi^2}{2c^2s^2}\biggr).\nonumber
\end{eqnarray}
Finally we have $c \leq\sqrt{s}/b$. Therefore, part (2) is proven
by combining (\ref{minimum}), (\ref{scadbdelta5})
(\ref{scadbdelta6}) and (\ref{scadbdelta7}), (\ref{scadbdelta8}).
This completes the proof.
\end{pf*}

\begin{pf*}{Proof of Theorem \ref{thm2}}
The proof is relegated to a supplementary file [\citet
{XueZouCai10}] for the sake of space.
\end{pf*}

\begin{pf*}{Proof of Corollary \ref{thm3}}
It follows directly from Theorems \ref{thm1}\break and~\ref{thm2}; thus we
omit its proof here.\vadjust{\goodbreak}
\end{pf*}

\begin{pf*}{Proof of Theorem \ref{thm4}}
Under the event $\Vert\widetilde\boldbeta{}^{(0)}-\boldbeta^*\Vert_{\infty
} \le\lambda$, we have $|\widetilde\boldbeta{}^{(0)}_{jk}|
\le\lambda$ for $(j,k) \in{\mathcal A^c}$
and $|\widetilde\boldbeta{}^{(0)}_{jk}| \ge a\lambda$ for $(j,k) \in
{\mathcal A}$. Therefore,\vspace*{1pt} $\widetilde\boldbeta{}^{(1)}$ is the solution
of the following penalized composite likelihood:
%
%
\begin{equation}
\label{thm4pfeq1} \widehat{\boldbeta}{}^{(1)}=\mathop{\arg\max}_{\boldbeta}
\biggl\{\ell_c(\boldbeta)-\lambda\sum_{(j,k)\in\mathcal A^c}|
\beta_{jk}|\biggr\}.
\end{equation}
It turns out that $\widehat{\boldbeta}{}^{\mathrm{oracle}}$ is the global
solution of (\ref{thm4pfeq1}) under the additional probability event
that $\{\Vert{\psi_{{\mathcal A}^c}(\widehat\boldbeta{}^{\mathrm{oracle}})}\Vert
_{\infty}
\leq\lambda\}$. To see this, we observe that for any $\boldbeta$,
\begin{eqnarray*}
&&\biggl(-\ell_c(\boldbeta)+\lambda\sum
_{(j,k)\in\mathcal A^c}|\beta_{jk}|\biggr)-\biggl(-\ell_c
\bigl(\widehat\boldbeta{}^{\mathrm{oracle}}\bigr)+\lambda\sum
_{(j,k)\in\mathcal A^c}\bigl|\widehat\beta{}^{\mathrm{oracle}}_{jk}\bigr|\biggr)
\\[-2pt]
&&\qquad\ge \sum_{(j,k) \in\mathcal{A}^c}\bigl(\lambda-\bigl|\psi_{(jk)}
\bigl(\widehat\boldbeta{}^{\mathrm{oracle}}\bigr)\bigr|\bigr) \cdot|\beta_{jk}|
\\[-2pt]
&&\qquad\ge 0,
\end{eqnarray*}
where we used the convexity of $-\ell_c$. In the proof of Theorem \ref
{thm1} we have shown that
\begin{eqnarray*}
&& \Pr\bigl(\bigl\Vert{ \psi_{{\mathcal A}^c}\bigl(\widehat\boldbeta{}^{\mathrm{oracle}}
\bigr)}\bigr\Vert_{\infty} > \lambda\bigr)
\\[-2pt]
&&\qquad< K^2\exp\biggl(-\frac{N\lambda^2}{32(2\phi+1)^2}\biggr)+{\exp\biggl(-
\frac
{N\lambda
}{ 3B(2\phi+1) s}\frac{b^2}{8^3}\biggr)}\\[-2pt]
&&\qquad\quad{}+K^2s\exp\biggl(-
\frac
{Nb^2}{2s^3}\biggr)
\nonumber
\\[-2pt]
&&\qquad\quad{}+2s^2\biggl[\exp\biggl(-\frac{b^2N}{8s^3}\biggr)+\exp\biggl(-
\frac{N}{s^2}\frac
{b^2}{2}\biggr)+\exp\biggl(-\frac{N}{s^2}
\frac{B^2}{8}\biggr)\biggr]
\\[-2pt]
&&\qquad \equiv \tau_3.
\nonumber
\end{eqnarray*}
Therefore, the LLA--CMA algorithm finds the oracle estimator with
probability at least $1-\tau_3-\Pr(\Vert\widetilde\boldbeta
{}^{(0)}-\boldbeta^* \Vert_{\infty}>\lambda)$.
This proves part (1).

If we further consider the event $\{\min_{(j,k) \in\mathcal
{A}} |{\widehat\beta{}^{\mathrm{oracle}}_{jk}}|>a\lambda\}$. Then $\widetilde
\boldbeta{}^{(2)}$ is the solution of the following penalized composite
likelihood
$
\max_{\boldbeta} \{\ell_c(\boldbeta)-\lambda\sum_{(j,k)\in
\mathcal A^c}|\beta_{jk}|\},
$
which implies that $\widetilde\boldbeta{}^{(2)}=\widetilde\boldbeta
{}^{(1)}$, and hence the LLA loop will stop.
From (\ref{minimum}) we have obtained a probability bound for the
event of $\{\min_{(j,k) \in\mathcal
{A}} |{\widehat\beta{}^{\mathrm{oracle}}_{jk}}| \le a\lambda\}$ as follows:
\begin{eqnarray*}
&&\Pr\Bigl(\min_{(j,k) \in\mathcal{A}} \bigl|\widetilde\beta{}^{\mathrm{hmle}}_{jk}\bigr| \le
a\lambda\Bigr)
\\[-2pt]
&&\qquad\le\exp\biggl(-R^{2}_*\frac{b^2}{8^3}\biggr)+2s^2
\exp\biggl(-\frac{N}{s^2}\frac
{b^2}{2}\biggr)+2s^2\exp
\biggl(-\frac{N}{s^2}\frac{B^2}{8}\biggr)
\\[-2pt]
&&\qquad\equiv \tau_4.\vadjust{\goodbreak}
\end{eqnarray*}
Then we have $\widetilde\boldbeta{}^{(m)}=\widetilde\boldbeta{}^{(1)}={
\widehat\boldbeta{}^{\mathrm{oracle}}}$ for $m=2,3,\ldots$ which means
the LLA--CMA algorithm converges after two LLA iteration and finds the
oracle estimator with
probability at least $1-\tau_3-\Pr(\Vert\widetilde\boldbeta
{}^{(0)}-\boldbeta^* \Vert_{\infty}>\lambda)-\tau_4$.
Note that $\tau_3+\tau_4=\tau_2$. This proves part (2).
\end{pf*}

\begin{pf*}{Proof of Corollary \ref{thm5}}
Part (1) follows directly from Theorem~\ref{thm4}. We only prove part
(2). With the chosen $\lambda^{\mathrm{lasso}}$, Theorem~\ref{thm2} shows
that with probability tending to one,\vspace*{2pt} $\widehat\boldbeta
{}^{\mathrm{lasso}}_{\mathcal{A}}=\widetilde\boldbeta_{\mathcal{A}} $,
$\widehat\boldbeta{}^{\mathrm{lasso}}_{\mathcal{A}^c}=0 $ and
$\pr(\Vert\widetilde\boldbeta_{\mathcal{A}}- \boldbeta^{*}_{\mathcal
{A}}\Vert_{2} \le16\lambda^{\mathrm{lasso}}\sqrt{s}/b)
\rightarrow0$. Note that
$16\lambda^{\mathrm{lasso}}\sqrt{s}/b < \lambda^{\mathrm{scad}}$ and $\Vert\widetilde
\boldbeta_{\mathcal{A}}- \boldbeta^{*}_{\mathcal{A}}\Vert_{\infty
} \le\Vert\widetilde\boldbeta_{\mathcal{A}}- \boldbeta^{*}_{\mathcal
{A}}\Vert_{2}$, we then conclude
$\tau_0=\pr(\Vert\widehat
\boldbeta{}^{\mathrm{lasso}}- \boldbeta^{*}\Vert_{\infty}
\le\lambda^{\mathrm{scad}}) \rightarrow0$.
\end{pf*}
\end{appendix}

\section*{Acknowledgments}

We thank the Editor, Associate Editor and
referees for their helpful comments.

\begin{supplement}
\stitle{Supplementary materials for ``Non-concave penalized composite
likelihood estimation of sparse Ising models''}
\slink[doi]{10.1214/12-AOS1017SUPP} 
\sdatatype{.pdf}
\sfilename{aos1017\_supp.pdf}
\sdescription{In this supplementary file, we provide a complete
theoretical analysis of the LASSO-penalized composite likelihood
estimator for sparse Ising models.}
\end{supplement}

%

\printaddresses

\end{document}